\documentclass[11pt, twoside, leqno]{article}

\usepackage{amssymb,latexsym}
\usepackage{enumerate}
\usepackage{color}
\usepackage{amsmath, amsthm, amsfonts, amssymb, mathrsfs}
\newtheorem{theorem}{Theorem}[section]
\newtheorem{corollary}[theorem]{Corollary}
\newtheorem{lemma}[theorem]{Lemma}
\newtheorem{proposition}[theorem]{Proposition}

\theoremstyle{definition}

\def\dint{\displaystyle\int}
\def\dsum{\displaystyle\sum}
\def\dsup{\displaystyle\sup}

\numberwithin{equation}{section}
\def\dfrac{\displaystyle\frac}

\numberwithin{equation}{section}

\numberwithin{equation}{section}

\begin{document}

\title{{\LARGE\bf     {Variational inequalities for the commutators of rough operators with BMO functions}
 \footnotetext{\small
{{\it MR(2000) Subject Classification}:\ } 42B20, 42B25}
\thanks {The research was supported by
NSF of China (Grant: 11471033, 11371057, 11571160, 11601396), Thousand Youth Talents Plan  of China (Grant: 429900018-101150(2016)),  Funds for Talents of China (Grant: 413100002), the Fundamental Research Funds for the Central Universities (FRF-BR-16-011A, 2014KJJCA10) and SRFDP of China   (Grant: 20130003110003).}}}

\date{ }
\maketitle

\begin{center}
{\bf Yanping Chen \footnote{\small {Corresponding author.\ }}}\\
Department of Applied Mathematics,  School of Mathematics and Physics,\\
 University of Science and Technology Beijing,\\
Beijing 100083, The People's Republic of China \\
 E-mail: {\it yanpingch@126.com}
\vskip 0.1cm

{\bf Yong Ding}\\
School of Mathematical Sciences, Beijing Normal University,\\
Laboratory of Mathematics and Complex Systems (BNU),
Ministry of Education,\\
Beijing 100875,  The People's Republic of China \\
E-mail: {\it dingy@bnu.edu.cn}\vskip 0.1cm

{\bf Guixiang Hong}\\
School of Mathematics and Statistics, Wuhan University,\\
Wuhan 430072,  The People's Republic of China \\
E-mail: {\it guixiang.hong@whu.edu.cn}\vskip 0.1cm and

{\bf Honghai Liu}\\

School of Mathematics and Information Science,\\
Henan Polytechnic University,
Jiaozuo, Henan, 454003,  China.\\
E-mail: {\it hhliu@hpu.edu.cn}

\end{center}

\begin{center}
\begin{minipage}{135mm}{\small
{\bf ABSTRACT\ }
In this paper, starting with a relatively simple observation that the variational estimates of the commutators of the standard Calder\'on-Zygmund operators with the BMO functions can be deduced from the weighted variational estimates of the standard Calder\'on-Zygmund operators themselves, we establish similar variational estimates for the commutators of the BMO functions with rough singular integrals which do not admit any weighted variational estimates. The proof involves many Littlewood-Paley type inequalities with commutators as well as Bony decomposition and related para-product estimates.}
\end{minipage}
\end{center}
\vspace{0.3cm}

 \section{Introduction}\label{sect1}
Motivated by the modulus of continuity of Brownian motion, L\'epingle \cite{Lep76} established the first variational inequality for general martingales (see \cite{PiXu88} for a simple proof). Bourgain \cite{Bou89} is the first one who exploited L\'epingle's result to obtain corresponding variational estimates for the Birkhoff ergodic averages along subsequences of natural numbers and then directly deduce pointwise convergence results without previous knowledge that pointwise convergence holds for a dense subclass of functions, which are not available in some ergodic models. In particular, Bourgain's work \cite{Bou89} has initiated a new research direction in ergodic theory and harmonic
analysis. In \cite{JKRW98}, \cite{JRW03}, \cite{JRW00}, \cite{CJRW2000}, \cite{CJRW2002}, Jones and his collaborators systematically studied variational inequalities for ergodic averages and truncated singular integrals of homogeneous type. Since then many other publications came to enrich the literature on this subject (cf. e.g. \cite{GiTo04}, \cite{LeXu2}, \cite{DMT12}, \cite{JSW08}, \cite{MaTo},\cite{OSTTW12}, \cite{Hon1}). Recently, several works on weighted as well as vector-valued variational inequalities in ergodic theory and harmonic analysis have also appeared (cf. e.g. \cite{MTX1},\cite{MTX2}, \cite{KZ}, \cite{HLP}, \cite{HoMa1}, \cite{HoMa2}); and several results on $\ell^p(\mathbb Z^d)$-estimates of $q$-variations for discrete operators of Radon type have also been established (cf. e.g. \cite{Kra14}, \cite{MiTr14}, \cite{MST15}, \cite{MTZ14}, \cite{Zor14}).

\bigskip

Most of the operators considered in the previous cited papers are of homogeneous type, and it is still unknown whether variational inequalities hold for all singular integrals of convolution type (while it is true when the kernel is smooth enough in \cite{MST15}), let alone for all standard Calder\'on-Zygmund operators. In our another paper \cite{CDHX}, we consider the variational inequality for Calder\'on commutators---commutators of singular integrals or pseudo differential operators with Lipschitz functions, which are typical example of Calder\'on-Zygmund operators of non-convolution type. In the present paper, we establish variational estimates for the commutators of rough operators with BMO functions, which finds its motivation in a simple observation that that the variational estimates of the commutators of the standard Calder\'on-Zygmund operators with the BMO functions can be deduced from the weighted variational estimates of the standard Calder\'on-Zygmund operators themselves.

To illustrate the observation, let us introduce some notations and recall some notions. Given a family of complex numbers $\mathfrak{a}=\{a_t: t\in\mathbb{R}\}$ and $\rho\ge1$, the $\rho$-variation norm of the family $\mathfrak{a}$ is defined by
\begin{equation}\label{q-ver number family}
\|\mathfrak{a}\|_{V_\rho}=\sup\big(\sum_{k\geq1}
|a_{t_k}-a_{t_{k-1}}|^\rho\big)^{\frac{1}{\rho}},
\end{equation}
where the supremum runs over all increasing sequences $\{t_k:k\geq0\}$. It is trivial that
\begin{equation}\label{number contr ineq}
\|\mathfrak{a}\|_{L^\infty(\mathbb{R})}:=\sup_{t\in\mathbb R}|a_t|
\le\|\mathfrak{a}\|_{V_\rho}+|a_{t_0}|\quad\text{for}\ \ \rho\ge1,
\end{equation}
for some fixed $t_0$.

Let $0<\rho<\infty$. Given a family of Lebesgue measurable functions $\mathcal F=\{F_t:t\in\mathbb{R}_+\}$ defined on $\mathbb{R}^n$, we define the strong $\rho$-variation function $V_\rho(\mathcal F)$ of the family $\mathcal F$ as
\begin{align*}
V_\rho(\mathcal F)(x)=\sup\|(F_{t_k}(x)-F_{t_{k-1}}(x))_{k\geq1}\|_{\ell^{\rho}}, a.e.\, x\in\mathbb R^n,
\end{align*}
where the supremum runs over all increasing sequences $\{t_k:k\geq0\}$. Suppose $\mathcal{T}=\{{T}_t\}_{t>0}$ is a family of operators  on $L^p(\Bbb R^n)\, (1\le p\le\infty)$. The strong $\rho$-variation operator is simply defined as
$$V_\rho \mathcal T(f)(x)=\|\{T_t(f)(x)\}_{t>0}\|_{V_\rho},\quad\forall f\in L^p(\mathbb{R}^n).$$ Thus the operator $V_\rho \mathcal T$ sends functions on $\Bbb R^n$ to nonnegative functions on $\Bbb R^n$.
It is easy to observe from the definition of $\rho$-variation norm that for any $x$ if $V_\rho \mathcal T(f)(x)<\infty$, then $\{T_t(f)(x)\}_{t>0}$ converges when $t\rightarrow0$ or $t\rightarrow\infty$. In particular, if $V_\rho \mathcal T(f)$ belongs to some function spaces such as $L^p(\mathbb{R}^n)$ or $L^{p,\infty}(\mathbb{R}^n)$, then the sequence converges almost everywhere without any additional condition. This is why mapping property of strong $\rho$-variation operator is so interesting in ergodic theory and harmonic analysis. Also, by \eqref{number contr ineq}, for any $f\in L^p(\Bbb R^n)$ and $x\in\Bbb R^n$, we have
 \begin{equation}\label{control of maxi opera}
T^\ast(f)(x)\le V_{\rho}(\mathcal{T}f)(x)\quad\text{for}\ \ \rho\ge1,
\end{equation}
where $T^\ast$ is the maximal operator defined by
$$T^\ast(f)(x):=\sup_{t>0}|{T}_t(f)(x)|.$$

For $b\in BMO(\mathbb R^n)$ and $u\in \Bbb
N$. If $T$ is a linear  operator on some measurable
function space, then the $u$-th order commutator formed by $b$ and $T$ is
defined by $T_{ b,u
}f(x):=T((b(x)-b(\cdot))^uf)(x)$. Simply, define $T_b:=T_{b,1}.$
In 1976, Coifman, Rochberg and Weiss [\ref{CRW}] obtained  a
characterization of $L^p$-boundedness of the commutators $R_{j;b}$
generated by the Reisz transforms $R_j\,(j=1,\cdots,n,)$ and a BMO
function $b$. As an application of this characterization, a
decomposition theorem of the real Hardy space is given in this
paper. Moreover, the authors in [\ref{CRW}] proved also that if
$\Omega\in {\rm Lip}(S^{n-1}),$
 then the commutator $T_{\Omega;b}$ for $T_\Omega$ and a BMO function $b$
  is bounded on $L^p$ for $1<p<\infty,$ which is defined by
$$T_{\Omega;b}f(x)=\text{p.v.}\dint_{
{\Bbb R}^n}\dfrac{\Omega(x-y)}{|x-y|^{n}}\ (b(x)-b(y))f(y)dy.$$
In the same paper, Coifman, Rochberg and Weiss [\ref{CRW}] outlined a different approach, which is less direct but shows the close relationship between the  weighted inequalities of the operator $T$ and the  weighted inequalities of the commutator $T_b.$
 In 1993,
Alvarez,  Bagby, Kurtz and  P\'{e}rez [\ref{ABKP}] developed the
idea of [\ref{CRW}], and established a generalized boundedness criterion for the commutators of linear operators. As it is well-known that the commutators have played an important role in harmonic
analysis and PDE, for example in the theory of non-divergent
elliptic equations with discontinuous coefficients (see \cite{C2},\,\cite{BC},\,\cite{CD},\,\cite{CFL},\,\cite{CFL2},\,\cite{FR},\,\,\cite{Ta3}). Moreover, there is also
an interesting connection between the nonlinear commutator,
considered by Rochberg and Weiss in [\ref{RW}], and Jacobian mapping
of vector functions. They have been applied in the study of the
nonlinear partial differential equations (see [\ref{CLMS},\,\ref{GI}]).

In order to state our main results, let us first recall  the definition and some properties of $A_p$ weight on $\mathbb R$. Let $w$ be a non-negative locally integrable function defined on $\mathbb R^n$. We say $w\in A_1$ if there is a constant $C>0$ such that $M(w)(x)\le Cw(x)$, where $M$ is the classical Hardy-Littlewood maximal operator defined by
$$Mf(x)=\sup_{r>0}\frac 1{r^n}\int_{|y|\le r}|f(x-y)|dy.$$
 Equivalently, $w\in A_1$ if and only if there is a constant $C>0$ such that for any cube $Q$
\begin{equation}\label{eq def A1}
\frac1{|Q|}\int_Qw(x)dx\le C\inf_{x\in Q}w(x).
\end{equation}
For $1<p<\infty$, we say that $w\in A_p$ if there exists a constant $C>0$ such that
\begin{equation}\label{eq def Ap}
\sup_Q\bigg(\frac1{|Q|}\int_Qw(x)dx\bigg)\bigg(\frac1{|Q|}\int_Qw(x)^{1-p'}dx\bigg)^{p-1}\le C.
\end{equation}
The smallest constant appearing in \eqref{eq def A1} or \eqref{eq def Ap} is denoted by $[w]_{A_p}$.

Now we state our first result as follows.
\begin{theorem}\label{thm:L}  Let $b\in BMO(\mathbb R^n)$. Let $\mathcal T=\{T_t:t\in\mathbb{R}_+\}$ be a family of linear operators and $\mathcal T_b=\{T_{t;b}:t\in\mathbb{R}_+\}$ be the family of commutators formed by the linear operators and $b$. Let $1<p<\infty$, $\rho\ge1$, $\tau\in {\Bbb R}$ and $1<s<\infty$. Let $w$ be a locally integrable function such that $w^\tau\in
A_s.$ Then  $V_\rho\mathcal T_b $ is bounded on $L^p(w)$, that is,
 $$\|V_\rho\mathcal T_b (f)\|_{L^p(w)}\le C\|b\|_{\ast}\|f\|_{L^p(w)}$$
provided that
$V_\rho\mathcal T $ is bounded on $L^p(w)$.
\end{theorem}

Now we present four consequences of Theorem \ref{thm:L}. Let $K$ be a kernel on $\Bbb R^n\times \Bbb R^n\setminus\{(x,x)\}:x\in \Bbb R^n\}$. We will suppose that $K$ satisify the following regularity conditions. There exist two constants $\delta>0$ and $C>0$ such that
 \begin{equation}\label{K0}
|K(x,y)|\le \frac{C}{|x-y|^n},\,\,\, \hbox{for}\,\,\,x\neq y;
 \end{equation}
 \begin{equation}\label{K1}
|K(x,y)-K(z,y)|\le \frac{C|x-z|^\delta}{|x-y|^{n+\delta}},\,\,\, \hbox{for}\,\,\,|x-y|>2|x-z|;
 \end{equation}
\begin{equation}\label{K2}
|K(y,x)-K(y,z)|\le \frac{C|x-z|^\delta}{|x-y|^{n+\delta}}, \,\,\,\hbox{for}\,\,\,|x-y|>2|x-z|.
 \end{equation}
Let $T$ be the standard Calder\'on-Zygmund operator associated to the kernel $K$: For a Schwartz function $f$,
 \begin{equation*}
Tf(x):=\lim_{\varepsilon\rightarrow0^+}T_{\varepsilon}f,
 \end{equation*}
where $T_\varepsilon$ be the truncated operator
 \begin{equation*}
T_\varepsilon f(x)=\int_{|x-y|>\varepsilon}K(x,y)f(y)\,dy; \end{equation*} for $b\in BMO(\mathbb R^n),$ define $T_{\varepsilon;b}$ to be the truncated commutator \begin{equation*}
T_{\varepsilon;b} f(x)=\int_{|x-y|>\varepsilon}K(x,y)(b(x)-b(y))f(y)\,dy. \end{equation*} Denote by $\mathcal T=\{T_\varepsilon\}_{\varepsilon>0}$ and $\mathcal T_b=\{T_{\varepsilon;b}\}_{\varepsilon>0}.$
 Let $K$ be kernel on $\Bbb R^n$  satisfying \eqref{K0}-\eqref{K2}, and let $2<\rho<\infty.$ In \cite{MTX2}, they showed that if the operator $V_\rho\mathcal T$ is of type $(p_0,p_0)$ for some $1<p_0<\infty,$ then for $1<p<\infty$ and $w\in A_p$, $
 V_\rho\mathcal T$  is bounded on $L^p(w).$ Thus apply Theorem \ref{thm:L}, we get

\begin{corollary}\label{cor:1}  Let $b\in BMO(\mathbb R^n)$ and $K$ be kernel on $\Bbb R^n$  satisfying \eqref{K0}-\eqref{K2}. Let $2<\rho<\infty.$
 If $V_\rho\mathcal T$ is of type $(p_0,p_0)$ for some $1<p_0<\infty,$
then for $1<p<\infty$ and $w\in A_p,$ there
exists a constant $C$ such that
\begin{equation*}
 \|V_\rho\mathcal T_b(f)\|_{L^p(w)}\le
C\|b\|_{\ast}\|f\|_{L^p(w)}.\end{equation*}
In particular, if $K$ is of convolution type satisfying cancellation condition
$$\int_{\partial B(0,t)}K(x)dx=0,\;\forall t>0$$
and either
{\rm (i)} $|\nabla K(x)|\leq C\frac{1}{|x|^{n+1}}$ or {\rm(ii)} $K(x)=\frac{\Omega(x)}{|x|^{n}}$ with $\Omega\in Lip_\delta(\mathbf S^n)$, then the associated $V_\rho\mathcal T_b$ if of type $L^p(w)$.
\end{corollary}

We refer the reader to \cite{MST15} (resp. \cite{CJRW2002}) for the result $V_\rho\mathcal T$ is $L^2$-bounded  when $K$ satisfies (i) (resp. (ii)). On the other hand, Corollary \ref{cor:1} implies the main result of \cite{LW} where $T$ is the Hilbert transform.

Let $\varphi:\Bbb R^n\rightarrow[0,+\infty)$ be a radially decreasing integrable function. Let $\varphi_t(x)=\frac{1}{t^n}\varphi(\frac{x}{t})$ and $\Phi(f)(x)=\{\varphi_t\ast f(x)\}_{t>0}$. In \cite{MTX2}, they also showed that for $2<\rho<\infty,$ the operator $V_\rho\Phi$ is bounded on $L^p(w)$ for $1<p<\infty$
 and $w\in A_p.$ Denote by $\Phi_b (f)=\{\varphi_t((b(x)-b(\cdot))f)\}_{t>0}$. Thus apply Theorem \ref{thm:L}, we have
\begin{corollary}\label{cor:2}  Let $b\in BMO(\mathbb R^n).$  Let $\varphi:\Bbb R^n\rightarrow[0,+\infty)$ be a radially decreasing integrable function. Let $\varphi_t(x)=\frac{1}{t^n}\varphi(\frac{x}{t})$ and $\Phi_b (f)(x)=\{\varphi_t((b(x)-b(\cdot))f)(x)\}_{t>0}$.
For  $2<\rho<\infty$, $1<p<\infty$ and $w\in A_p,$ there
exists a constant $C$ such that
\begin{equation}\label{Vqr}
 \|V_\rho\Phi_b (f)\|_{L^p(w)}\le
C\|b\|_{\ast}\|f\|_{L^p(w)}.\end{equation}

\end{corollary}

\bigskip

Not only to smooth singular kernels or good approximation identities, Theorem \ref{thm:L} can also be applied  to singular integral operators or averaging operators with some homogeneous rough kernels.
Suppose $T_{ \Omega,\varepsilon}$ is the truncated singular integral operator defined by
 \begin{equation}\label{tr of S}
T_{ \Omega,\varepsilon}
f(x)=\int_{|y|>\varepsilon}\frac{\Omega(y')}{|y|^n}f(x-y)dy,
 \end{equation}
where $\Omega\in L^1({\mathbf S}^{n-1})$
 satisfies the cancelation condition
\begin{equation}\label{can of O}
\int_{\mathbf S^{n-1}}\Omega(y')d\sigma(y')=0.
 \end{equation}
 For $1<p<\infty$ and $f\in C^\infty_c(\mathbb R^n)$,
 the Calder\'on-Zygmund singular integral operator $T$ with homogeneous kernel is defined by
\begin{equation}\label{SIO}
T_{ \Omega}f(x)=\lim_{\varepsilon\rightarrow0^+}T_\varepsilon f(x), \ a.e. \ x\in\mathbb R^n.
\end{equation}
Denote the family of operators $\{T_{ \Omega,\varepsilon}\}_{\varepsilon>0}$  by $\mathcal T_\Omega$.
 In \cite{CDHL}, we showed that if $\Omega\in L^q(\mathbf S^{n-1})$, $q>1$ satisfying \eqref{can of O},
then for $\rho>2$,
$V_\rho\mathcal T_\Omega$ is bounded on $L^p(w),$ whenever $w$ and $p$ satisfy one of the following conditions:
\begin{enumerate}[(i)]
\item  $q'\le p<\infty,\, p\neq 1$  and $ w\in A_{p/q'},$
\item $1< p\le q,\, p\neq \infty$  and $ w^{-\frac{1}{(p-1)}}\in A_{p'/q'}$.
\end{enumerate}
Thus, we can invoke Theorem \ref{thm:L} to obtain the corresponding weighted estimates for the strong $\rho$ variation of commutators with rough kernel and $b\in BMO(\mathbb R^n). $

\begin{corollary}\label{cor:3}  Let $b\in BMO(\mathbb R^n)$ and $\mathcal T_{\Omega;b}f(x)$ be the family of commutators $\{T_{ \Omega,\varepsilon}((b(x)-b(\cdot)f)(x))\}_{\varepsilon>0}$ with  $\Omega\in L^q(\mathbf S^{n-1})$, $q>1$ satisfying \eqref{can of O},
then for $\rho>2$, there
exists a constant $C$ such that
\begin{equation}\label{Vqr}
 \|V_\rho\mathcal T_b(f)\|_{L^p(w)}\le
C\|b\|_{\ast}\|f\|_{L^p(w)},\end{equation}
if $w$ and $p$ satisfy one of the situation (i) or (ii).

\end{corollary}

\bigskip

We can also give  applications of Theorem \ref{thm:L} to the situation of averaging operators with rough kernels $\mathcal M_\Omega=\{M_{\Omega,t}\}_{t>0}$, where $M_{\Omega,t}$ is defined as
\begin{equation}\label{def M O}
M_{\Omega,t} f(x)=\frac1{t^n}\int_{|y|<t}\Omega(y')f(x-y)dy,
\end{equation}
where $\Omega\in L^1({\mathbf S}^{n-1})$.
In \cite{CDHL}, we showed that if $\Omega\in L^q(\mathbf S^{n-1})$, $q>1$,
then for $\rho>2$,
$
V_\rho\mathcal M_\Omega $ is bounded on $L^p(w)$
if $w$ and $p$ satisfy one of the  conditions (i) and (ii). See
\cite{GHST,CD} on the maximal inequality for the family $\mathcal M_{\Omega}$.
Given a BMO function $b$, denote the family of operators $\{M_{\Omega,t}(b(x)-b(\cdot)f)(x)\}_{t>0}$  by  $\mathcal M_{\Omega;b}f(x)$. Thus, we can invoke Theorem \ref{thm:L} to obtain the corresponding weighted estimates for the strong $\rho$ variation of commutators with rough kernel and $b\in BMO(\mathbb R^n). $

\begin{corollary}\label{cor:4}  Let $b\in BMO(\mathbb R^n)$ and $\mathcal M_{\Omega;b}$ be the family of the commutators of averaging operators with  $\Omega\in L^q(\mathbf S^{n-1})$, $q>1$.
Then for $\rho>2$, there
exists a constant $C$ such that
\begin{equation}\label{Vqr}
 \|V_\rho\mathcal M_{\Omega;b}(f)\|_{L^p(w)}\le
C\|b\|_{\ast}\|f\|_{L^p(w)},\end{equation}
if $w$ and $p$ satisfy one of the  conditions (i) or (ii).

\end{corollary}

\bigskip

However,  it is not clear up to now whether the  operator $V_\rho\mathcal T_\Omega$ or $V_\rho\mathcal M_\Omega$
with $\Omega\in L^1\setminus\bigcup_{q>1}L^q(S^{n-1})$ is bounded on
$L^p(w)$ ($1<p<\infty$) for all $w\in A_r$ ($1<r<\infty$), Hence, if $\Omega\in L^1\setminus\bigcup_{q>1}L^q(S^{n-1}),$
the $L^p$ boundedness of $V_\rho\mathcal T_{\Omega;b}$ and $V_\rho\mathcal M_{\Omega;b}$ can not be deduced from
Theorem \ref{thm:L} .

The main purpose of this paper is to give a sufficient condition which contains $\bigcup_{q>1}L^q(S^{n-1})$, such that the operators $V_\rho\mathcal T_{\Omega,b}$  and $V_\rho\mathcal M_{\Omega;b}$ are bounded on $L^p({\Bbb R}^n)$ for $1<p<\infty.$ It is well known that $$\bigcup_{q>1}L^q(S^{n-1})\subset L(\log^+L)^\alpha(S^{n-1})$$ for any $\alpha>0.$

The second main result of this paper is formulated as follows.

\begin{theorem}\label{thm:M}  Let $b\in BMO(\mathbb R^n)$ and $\mathcal T_{\Omega;b}$ be the family of the commutators of truncated singular integral operators with $\Omega$ satisfying \eqref{can of O}.
 If $\Omega\in L(\log^+\!\! L)^3(\mathbf S^{n-1})$, then the following  $\rho\,(2<\rho<\infty)$-variational inequality holds for $1< p<\infty,$ namely,
$$
\|V_{\rho}\mathcal T_{\Omega;b} (f)\|_{L^p}\le
C\|b\|_{\ast}\|f\|_{L^p}.$$
\end{theorem}

The proof of this result is based on Fourier transform, which is somehow standard but technical since it involves many Littlewood-Paley type inequalities with commutators as well as Bony decomposition and related para-product estimates.

Our approach to the variational estimates for singular integrals also works for the family $\mathcal M_{\Omega;b}$.
\begin{theorem}\label{thm:N}
Let $1<p<\infty$ and  Let $b\in BMO(\mathbb R^n)$. If $\Omega\in L(\log^+\!\!L)^{2}(\mathbf S^{n-1})$,
then the  $\rho\,(2<\rho<\infty)$-variational inequality for the family $\mathcal M_{\Omega;b}$ holds,
$$
\|V_{\rho}\mathcal M_{\Omega;b} (f)\|_{L^p}\le
C\|b\|_{\ast}\|f\|_{L^p},\;\forall f\in L^p(\mathbb R^n).
$$
\end{theorem}

The paper is organized as follows. In Section 2,   we give the proof of Theorem \ref{thm:L}.
In Section 3,  we give the proof of Theorem  \ref{thm:M}.   Section 4 is devoted to  the proof of  Theorem \ref{thm:N}.
For $p\ge 1,$ $p'$ denotes the conjugate
exponent of $p$, that is, $p'=p/(p-1).$  Throughout this paper, the letter $``C\,"$ will stand for a positive
constant which is independent of the essential variables and not
necessarily the same one in each occurrence.

\section{Proof of Theorem \ref{thm:L}}
Denote $F(z)=e^{z[b(x)-b(y)]},\,z\in {\Bbb C}.$  Then by the analyticity of $F(z)$ on $\Bbb C$ and the Cauchy integration formula, we have for any $\varepsilon>0$
\begin{align*}
b(x)-b(y)=F'(0)&=\dfrac{1}{2\pi i}\int_{|z|=\varepsilon}\dfrac{F(z)}{z^2}\,dz\\&=\dfrac{1}{2\pi }\dint_{0}^{2\pi}\dfrac{F(\varepsilon e^{i\theta})}{\varepsilon e^{i\theta}}\,d\theta\\&=\dfrac{1}{2\pi\varepsilon  }\dint_{0}^{2\pi}{e^{\varepsilon e^{i\theta}[b(x)-b(y)]}} e^{-i\theta}\,d\theta.
\end{align*}
For any linear operator $T$ we get
\begin{align*}
T_b f(x)=T(b(x)-b(\cdot)f)(x)&=\dfrac{1}{2\pi\varepsilon  }\dint_{0}^{2\pi}T(e^{-b(\cdot)\varepsilon e^{i\theta}}f)(x)e^{b(x)\varepsilon e^{i\theta}}e^{-i\theta }\,d\theta\\&=\dfrac{1}{2\pi\varepsilon  }\dint_{0}^{2\pi}T(h_\theta)(x)e^{b(x)\varepsilon e^{i\theta}}e^{-i\theta }\,d\theta
\end{align*}
where $h_\theta(x)=f(x)e^{-b(x)\varepsilon e^{i\theta}}$ for $\theta\in [0,2\pi].$
Recall that
\begin{equation*}
V_\rho(\mathcal F)(x)=\sup\|(F_{t_k}(x)-F_{t_{k-1}}(x))_{k\geq1}\|_{\ell^{\rho}},\end{equation*}where the supremum runs over all increasing sequences $\{t_k:k\geq0\}$.
Then using the Minkowski inequality, we have for $w^\tau\in A_{s}$
\begin{align}\label{Vqb}
&\|V_\rho\mathcal T_b (f)\|_{L^p(w)}\\&\nonumber=\|\sup\|(T_{t_k;b}f-T_{t_{k-1};b}f)_{k\geq1}\|_{\ell^{\rho}}\|_{L^p(w)}\\&
\nonumber=\bigg(\dint_{{\Bbb R}^n}\bigg|\sup\Big\|\big(\dfrac{1}{2\pi\varepsilon  }\dint_{0}^{2\pi}\big(T_{t_k}(h_\theta)(x)-T_{t_{k-1}}(h_\theta)(x)\big)e^{b(x)\varepsilon e^{i\theta}}e^{-i\theta }\,d\theta\big)_{k\geq1} \Big\|_{\ell^\rho}\bigg|^{p}w(x)\,dx\bigg)^{1/p}\\&\nonumber
\le \dfrac{1}{2\pi\varepsilon  }\dint_{0}^{2\pi}\bigg(\dint_{{\Bbb R}^n}\bigg|\sup\|\big(T_{t_k}(h_\theta)(x)-T_{t_{k-1}}(h_\theta)(x)\big)_{k\geq1}\|_{\ell^\rho}\bigg|^{p}e^{pb(x)\varepsilon cos\theta}w(x)\,dx\bigg)^{1/p} \,d\theta\\&
\nonumber= \dfrac{1}{2\pi\varepsilon  }\dint_{0}^{2\pi}\|V_\rho\mathcal T( h_\theta)\|_{L^p(we^{pb\varepsilon cos\theta})}\,d\theta.
\end{align}
 Note that for $f\in L^p(w),$ it is easy to check that for any $\theta\in [0,2\pi],$
$$h_\theta\in L^p(we^{pb(\cdot)\varepsilon \cos\theta})\,\,\, \hbox{and}\,\,\,\,\|h_\theta\|_{L^p(we^{pb(\cdot)\varepsilon \cos\theta})}=\|f\|_{L^p(w)}.$$ Hence we should compute
\begin{align}\label{Vqw}[e^{\tau pb\varepsilon \cos
\theta}w^\tau]_{A_s}&=\sup\limits_{Q}\bigg(\dfrac{1}{|Q|}\dint_{Q}
e^{\tau pb(x)\varepsilon \cos \theta}w^\tau(x)\,dx\bigg)\bigg(\frac{1}{|Q|}\dint_Q
\bigg[e^{\tau pb(x)\varepsilon \cos
\theta}w^\tau(x)\bigg]^{-\frac{1}{s-1}}\,dx\bigg)^{s-1}.\end{align}
Now, since $w^\tau \in A_s$,
there exists some $r_1>1$ such that
$$\bigg(\frac{1}{|Q|}\dint_{Q}w^{\tau r_1}\,dx\bigg)^{\frac{1}{r_1}}\leq \frac{2}{|Q|}\dint_{Q}w^\tau\,dx$$
and similarly for $w^{-\frac{\tau}{s-1}} \in A_{s'},$ there exists some $r_2>1$ such that
 then
$$\bigg(\frac{1}{|Q|}\dint_{Q}w^{-\frac{\tau{r_2}}{s-1}}\,dx\bigg)^{\frac{1}{r_2}}\leq \frac{2}{|Q|}\dint_{Q}w^{-\tau/{(s-1)}}\,dx.$$ By this, we know if $r_2<r_1,$ by the H\"{o}lder inequality, we get
\begin{align}\label{Vqr}\bigg(\frac{1}{|Q|}\dint_{Q}w^{\tau r_2}\,dx\bigg)^{\frac{1}{r_2}}\le \bigg(\frac{1}{|Q|}\dint_{Q}w^{\tau r_1}\,dx\bigg)^{\frac{1}{r_1}}
\le \frac{2}{|Q|}\dint_{Q}w^\tau\,dx.\end{align} If
$r_1<r_2,$ by the H\"{o}lder inequality, we get
\begin{align}\label{Vqr1}\bigg(\dfrac{1}{|Q|}\dint_{Q}w^{-\frac{\tau r_1}{s-1}}\,dx\bigg)^{\frac{1}{r_1}}\le \bigg(\frac{1}{|Q|}\dint_{Q}w^{-\frac{\tau r_2}{s-1}}\,dx\bigg)^{\frac{1}{r_2}}
\le \frac{2}{|Q|}\dint_{Q}w^{-\frac{\tau}{s-1}}\,dx.\end{align} Taking
$r=\min\{r_1,r_2\},$ using \eqref{Vqr}, \eqref{Vqr1} and the H\"{o}lder inequality we have
 \begin{align}\label{Vqwe}
&[e^{\tau pb\varepsilon \cos
\theta}w^\tau]_{A_s}\\&\nonumber=\sup_{Q}\bigg(\dfrac{1}{|Q|}\int_{Q}e^{\tau pb(x)\varepsilon
\cos\theta}w^\tau(x)\,dx\bigg)\bigg(\dfrac{1}{|Q|}\int_{Q}
e^{-\tau\frac{s}{s-1}b(x)\varepsilon
\cos\theta}w(x)^{-\frac{\tau}{s-1}}\,dx\bigg)^{s-1}\\&\le\nonumber\sup_{Q}\bigg(\dfrac{1}{|Q|}\int_{Q}
w^{\tau r}(x)\,dx\bigg)^{\frac{1}{r}}\bigg(\dfrac{1}{|Q|}\int_{Q}e^{\tau pr'b(x)\varepsilon
\cos\theta}\,dx\bigg)^{\frac{1}{r'}}\\&\nonumber\quad\bigg(\dfrac{1}{|Q|}\int_{Q}
w(x)^{-\tau r
/(s-1)}\,dx\bigg)^{\frac{s-1}{r}}\bigg(\dfrac{1}{|Q|}\int_{Q}
e^{-\frac{\tau p}{s-1}r'b(x)\varepsilon
\cos\theta}\,dx\bigg)^{\frac{s-1}{r'}}\\&\nonumber\le
4\sup_{Q}\bigg(\dfrac{1}{|Q|}\int_{Q}
w^\tau(x)\,dx\bigg)\bigg(\dfrac{1}{|Q|}\int_{Q} w(x)^{-\tau
/(s-1)}\,dx\bigg)^{s-1}\\&\nonumber\sup_{Q}\bigg(\dfrac{1}{|Q|}\int_{Q}e^{\tau pr'b(x)\varepsilon
\cos\theta}\,dx\bigg)^{\frac{1}{r'}}\bigg(\dfrac{1}{|Q|}\int_{Q}
e^{-\frac{\tau p}{s-1}r'b(x)\varepsilon
\cos\theta}\,dx\bigg)^{\frac{s-1}{r'}}\\&\nonumber=4[w^\tau]_{A_s}[e^{\tau pr'b\varepsilon
\cos\theta}]_{A_s}^{\frac{1}{r'}}.\end{align}
Now, since
$b\in BMO(\mathbb R^n)$ we  apply Lemma 2.1 to do this.

\begin{lemma}\label{esb}(\cite{CPP})\ Let $1<p<\infty$ and  $b\in BMO$. There exist two dimensional constants $\alpha_n, \beta_n$ satisfying $0<\alpha_n<1$ and $0<\beta_n<\infty$ such that
for any $\lambda\in \Bbb R$ with $|\lambda|\le \frac{\alpha_n}{\|b\|_{\ast}}\min\{1,\frac{1}{p-1}\}$, we have
$$e^{\lambda b}\in A_p\,\,\,\hbox{and}\,\,\,[e^{\lambda b}]_{A_p}\le \beta_n^p.$$
\end{lemma}

Now, we choose the radius
\begin{align}\label{ese}\varepsilon=\frac{\alpha_n\min\{1,\frac{s'}{s}\}}{\tau pr'\|b\|_{\ast}},\end{align} such that $|\tau pr'\varepsilon \cos\theta|\le \frac{\alpha_n}{\|b\|_{\ast}}\min\{1,\frac{s'}{s}\}$, then by Lemma \ref{esb}, we get
\begin{align}\label{esp}[e^{\tau pr'b\varepsilon \cos\theta}]_{A_s}\le \beta_n^s.\end{align}
 Combining \eqref{Vqwe} and \eqref{esp},  we get
\begin{align}\label{esc}
[e^{\tau pb\varepsilon \cos\theta}w^\tau]_{A_s}\le 4[w^\tau]_{A_s}[e^{\tau pr'b\varepsilon \cos\theta}]_{A_s}^{\frac{1}{r'}} \le 4[w^\tau]_{A_s}\beta_n^{s/r'} \le C_1[w^\tau]_{A_s},\end{align} where $C_1$ is independent of $w$ and $b.$
Now we return to \eqref{Vqb}. By \eqref{ese} and \eqref{esc}, we get for $1<p<\infty$ and $w^\tau\in A_{s}$,
\begin{align*}
\|V_\rho\mathcal T_b (f)\|_{L^p(w)}&\leq
\dfrac{C}{2\pi\varepsilon}\dint_0^{2\pi}\|h_\theta\|_{L^p(e^{pb(\cdot)\varepsilon
\cos \theta}w)}\,d\theta\leq C\|b\|_{\ast} \|f\|_{L^p(w)}.
\end{align*}
Thus we complete the proof of Theorem \ref{thm:L}.
\qed

\section{Proof of Theorem \ref{thm:M}}

As having appeared in the proof of most of variational inequalities (in particular see \cite{JSW08}), we shall show the desired  estimate by proving separately the long and short variational estimates. That is, we are reduced to show
\begin{align}\label{dyadic jump singular}
\|V_{\rho}(\{T_{\Omega,2^k;b}f\}_k)\|_{L^p}\le
C\|b\|_{\ast}\|f\|_{L^p}
\end{align}
and
\begin{align}\label{short variation singular}
\|S_2(\mathcal T_{\Omega;b}f)\|_{L^p}\le
C\|b\|_{\ast}\|f\|_{L^p},
\end{align}
where
 $$
S_2(\mathcal T_{\Omega;b}f)(x)=\bigg(\sum_{j\in\mathbb Z}[V_{2,j}(\mathcal T_{\Omega;b} f)(x)]^2\bigg)^{1/2},
$$
with
$$
V_{2,j}(\mathcal T_{\Omega;b}f)(x)=\bigg(\sup_
{2^j\leq t_0<\cdots<t_N<2^{j+1}}\sum_{l=0}^{N}| T_{\Omega,t_{k+1};b}f(x)- T_{\Omega,t_k;b}f(x)|^2\bigg)^{1/2}.
 $$

To deal with the long variation \eqref{dyadic jump singular}, we use weighted estimates obtained in our previous paper \cite{CDHL} together with some Littlewood-Paley type estimates involving commutators, in addition to Theorem \ref{thm:L}.
While to handle with the short variation \eqref{short variation singular}, we exploit Bony decomposition and paraproduct estimates, as well as Fefferman-Stein inequality for rough maximal functions.

\subsection{Proof of Theorem \eqref{dyadic jump singular}}

In this subsection, we give the proof of \eqref{dyadic jump singular} and in the next subsection we deal with \eqref{short variation singular}. Let us begin with one definition. For
$j\in\mathbb Z$,  let $\nu_j(x)=\frac{\Omega(y)}{|y|^n}\chi_{\{2^j\le |x|<2^{j+1}\}}(x)$, then
$$
\nu_j\ast f(x)=\int_{2^j\le |y|<2^{j+1}}\frac{\Omega(y)}{|y|^n}f(x-y)dy.
$$
Let $\phi\in \mathscr{S}(\mathbb R^n)$ be a radial function such that $\hat{\phi}(\xi)=1$ for
$|\xi|\le2$ and $\hat{\phi}(\xi)=0$ for $|\xi|>4$. We have the following decomposition
\begin{align*}
T_{\Omega,2^k}f&=\phi_k\ast T_\Omega f+\sum_{s\ge0}(\delta_0-\phi_k)\ast\nu_{k+s}\ast f-\phi_k\ast\sum_{s<0}\nu_{k+s}\ast
f,
\end{align*} where $\phi_k$ satisfies $\widehat{\phi_k}(\xi)=\hat{\phi}(2^k\xi)$, $\delta_0$ is the Dirac measure at 0 and $s\in\Bbb N\cup \{0\}$. Then we have \begin{align*}
T_{\Omega,2^k;b}f&=(\phi_k\ast T_\Omega)_b f+\Big(\sum_{s\ge0}(\delta_0-\phi_k)\ast\nu_{k+s}\Big)_b f-\Big(\phi_k\ast\sum_{s<0}\nu_{k+s}\Big)_b
f\\
&:=T^1_{k;b}f+T^2_{k;b}f-T_{k;b}^3f.
\end{align*}
 Let $\mathscr T^i_bf$ denote the family $\{T^i_{k;b}f\}_{k\in\mathbb Z}$  for
$i=1,2,3$. Obviously, to show \eqref{dyadic jump singular} it suffices to prove the following inequalities:
\begin{equation}\label{lnti}
\|V_\rho(\mathscr T_b^if)\|_{L^p}\le C\|b\|_{\ast}\|f\|_{L^p},\ \ 1<p<\infty,\ \ i=1,2,3.
\end{equation}

\bigskip
\noindent {\bf Estimate of \eqref{lnti} for $i=1$.}
We need the following two lemmas.
\begin{lemma}\label{Phi} (\cite{DHL})  Let $\phi_k$ be given above, set $\mathscr Uf=\{\phi_k\ast f\}_{k}$. Then for $1<p<\infty$ and $2<\rho<\infty$,
$$
\|V_\rho(\mathscr Uf)\|_{L^p}\le C\|f\|_{L^p}.
$$
\end{lemma}
In \cite{CDHL}, we proved that for $2<\rho<\infty,$ $1<p<\infty$ and $w\in A_p,$$
V_\rho\mathscr U$ is bounded on $L^p(w).$ Then apply Theorem \ref{thm:L}, we get
\begin{lemma}\label{Phibj}
For $k\in \Bbb Z,$ let $\phi_k$ be given above and $b\in BMO({\Bbb R}^n)$. Denote by $\Phi_kf=\phi_k\ast f,\, k\in \Bbb Z.$ Let $\mathscr U_bf=\{\Phi_{k;b} f\}_{k}$. Then for $1<p<\infty$ and $w\in A_p$,
$$
\|V_\rho\mathscr U_b(f)\|_{L^p(w)}\le C\|f\|_{L^p(w)}.$$
\end{lemma}
 Write
 $$T^1_{k;b}f=\Phi_{k;b} T_\Omega f+\Phi_k T_{\Omega;b}f.$$
Then combining Lemma \ref{Phi} with Lemma \ref{Phibj}, the $L^p$-boundedness of $T_\Omega$ (see \cite{CZ}) and $[b,T_\Omega]$ (see \cite{H}), we can get the following estimate easily
\begin{align*}
\|V_q(\mathscr T_b^1f)\|_{L^p}&\le\|V_q(\{\Phi_k (T_{\Omega;b}f)\}_k)\|_{L^p}+\|V_q(\{\Phi_{k;b}  (T_\Omega f)\}_k)\|_{L^p}\\&\le C\bigg(\| T_{\Omega;b}f\|_{L^p}+\|b\|_{\ast}\| T_\Omega f\|_{L^p}\bigg)\\&\le C\|\Omega\|_{L({\log^+L})^2(\mathbf S^{n-1})}\|b\|_{\ast}\|f\|_{L^p}.
\end{align*}

\medskip

\noindent {\bf Estimate of \eqref{lnti} for $i=2$.} We first give two lemmas will be used in the following arguments.
\begin{lemma}\label{mulcom}(\cite{H}) Let $m\in C_0^\infty({\Bbb
R}^n)$ and
$\hbox{supp}\ m\subset \{ |\xi|\le
2\sigma\}$ for some $\sigma\in(0,\infty)$. Suppose $m$ satisfies
$$\|m\|_{L^\infty}\le C 2^{ -\gamma s}\min\{\sigma,\sigma^{-\lambda}\},\ \ \
\|\nabla m\|_{L^\infty}\le C2^s$$
for some constants $C$, $\lambda,\gamma>0$ and $s\in \Bbb N$.
 Let $T_m$ be the multiplier operator defined by
$\widehat{T_m f}(\xi)=m(\xi)\widehat{f}(\xi)$. Moreover, for $b\in BMO$ and $u\in \Bbb
N,$ denote by $T_{m; b,u
}f(x)=T_{m}((b(x)-b(\cdot))^uf)(x)$ the $u$-th order
commutator of $T_{m}.$
Then for  any $0<v<1$,
there exist positive constants $C=C(n,v)$ and  $\beta\in (0,1)$ such that
$$\|T_{m;b,u}f\|_{L^{2}}\le C2^{-\beta s}\min\{\sigma^v,\sigma^{-\lambda v}\} \|b\|_{\ast}^u\|f\|_{L^{2}}.$$
\end{lemma}\begin{lemma}\label{Deltacom} (see \cite{H}). Let $b\in BMO(\mathbb R^n)$. Then, for  $1<p<\infty$ and
$f\in L^p({\Bbb R}^n)$, we have

\rm(\rm i) \quad $\bigg\|\bigg(\dsum_{l\in \Bbb
Z}|\Delta_{l;b,u}f|^2\bigg)^{1/2}\bigg\|_{L^p}\le
C(n,p)\|b\|_{\ast}^u \|f\|_{L^p};$

\rm(\rm ii)\quad $\bigg\|\bigg(\dsum_{l\in \Bbb
Z}|\Delta_{l;b,u}^2f|^2\bigg)^{1/2}\bigg\|_{L^p}\le
C(n,p)\|b\|_{\ast}^u \|f\|_{L^p}$.
\end{lemma}

Let $E_0=\{ x^\prime\in \mathbf S^{n-1}:
|\Omega(x^\prime)|<2\}$
 and $E_{d}=\{ x^\prime\in \mathbf S^{n-1}:
 2^{d}\le|\Omega(x^\prime)|<2^{d+1}\}$
  for positive integer $d$.  For $d\ge 0,$ let
$$\Omega_d(y')=\Omega(y')\chi_{E_d}(y')-\frac{1}
{|\mathbf S^{n-1}|}\int_{E_d}\Omega(x' ) \,d\sigma(x').$$ Since $\Omega$
satisfies \eqref{can of O}, then
$$\int_{\mathbf S^{n-1}}\Omega_d(y')\,d\sigma(y')=0 \
 \ \hbox{for}\ \  d\ge 0$$
and $\Omega(y')=\sum_{d\ge 0}\Omega_d(y').$ Set
   $$\nu_{j,d}(x)=\frac{\Omega_d(x)}
   {|x|^{n}}\chi_{\{2^j\le |x|
   < 2^{j+1}\}}(x).$$
Then, by the Minkowski inequality, we get
\begin{align}
\nonumber V_q\mathscr T^{2}_b(f)(x)&\le\dsum_{s\ge0}\Big(\dsum_{k\in \Bbb Z}\Big|\Big((\delta_0-\phi_k)\ast\nu_{k+s}\Big)_b
f(x)\Big|^2\Big)^{1/2}\\&\le\dsum_{s\ge0}\dsum_{d\ge0}\Big(\dsum_{k\in \Bbb Z}\Big|\Big((\delta_0-\phi_k)\ast\nu_{k+s,d}\Big)_b
f(x)\Big|^2\Big)^{1/2}.
\end{align}
Let $\varphi\in
C_0^\infty({\Bbb R}^n)$ be a radial function such
  that
 $0\le \varphi\le 1,$ $\hbox{supp}\ \varphi \subset\{1/2\le |\xi|\le
2\}$ and $\sum_{l\in \Bbb Z}\varphi^2(2^{-l}\xi)=1$ for $|\xi|\neq 0.$
Define the multiplier $\Delta_l$ by
$ \widehat{\Delta_l f}(\xi)=\varphi(2^{-l}\xi)\widehat{f}(\xi).$
By the Minkowski inequality again, we get
\begin{align}
\nonumber V_\rho\mathscr T^{2}_b(f)(x)&\le\dsum_{l\in \Bbb Z}\dsum_{s\ge0}\dsum_{d\ge0}\Big(\dsum_{k\in \Bbb Z}\Big|\Big((\delta_0-\phi_k)\ast\nu_{k+s,d}\ast\Delta_{l-k}^2\Big)_b
f(x)\Big|^2\Big)^{1/2}.
\end{align}
We set
$$m_{s,k,d}(\xi)=(1-\widehat{\phi_k}(\xi))\widehat{\nu_{k+s,d}}(\xi),\,\,\,\,\,m_{s,k,d}^l(\xi)=(1-\widehat{\phi_k}(\xi))\widehat{\nu_{k+s,d}}(\xi)\varphi(2^{k-l}\xi)$$
Define multipliers $F_{s,k,d}$  and $F_{s,k,d}^l$, respectively by $$\widehat{F_{s,k,d}f}(\xi)=m_{s,k,d}(\xi)\widehat{f}(\xi)$$ and $$\widehat{F_{s,k,d}^lf}(\xi)=m_{s,k,d}^l(\xi)\widehat{f}(\xi).$$
Then \begin{align}
\nonumber V_\rho\mathscr T^{2}_b(f)(x)&\le\dsum_{l\in \Bbb Z}\dsum_{s\ge0}\dsum_{d\ge0}\Big(\dsum_{k\in \Bbb Z}\Big|\big( F_{s,k,d}\Delta_{l-k}^2\big)_b
f(x)\Big|^2\Big)^{1/2}.
\end{align}
If we can  prove that
\begin{align}\label{Fsk}
\bigg\|\bigg(\dsum_{k\in \Bbb Z}
 |\big( F_{s,k,d}\Delta_{l-k}^2\big)_b f|^2\bigg)^{1/2}\bigg\|_{L^2}
 &\le C\|\Omega_d\|_{L^\infty(\mathbf S^{n-1})}2^{-\beta s}2^{-\theta |l|}\|b\|_{\ast}\|f\|_{L^2}
\end{align} and for $1<p<\infty$
\begin{align}\label{Fsk1}\bigg\|\bigg(\dsum_{k\in \Bbb Z}
 |\big( F_{s,k,d}\Delta_{l-k}^2\big)_b f|^2\bigg)^{1/2}\bigg\|_{L^p}\le
 C\|\Omega_d\|_{L{\log^+L}(\mathbf S^{n-1})}\|b\|_{\ast}\|f\|_{L^p},
\end{align} then we may finish the  estimate of \eqref{lnti} for $i=2$.
In fact, interpolating \eqref{Fsk} and \eqref{Fsk1},  we get for $0<\theta_0, \beta_0<1$,
\begin{equation}\label{Fsk2}
\bigg\|\bigg(\dsum_{k\in \Bbb Z}
 |\big( F_{s,k,d}\Delta_{l-k}^2\big)_bf|^2\bigg)^{1/2}\bigg\|_{L^p}\le C 2^{-\beta_0 s}2^{- \theta_0|l|}
\|b\|_{\ast}\|\Omega_d\|_{L^\infty(\mathbf S^{n-1})}\|f\|_{L^p},\,\,\,1<p<\infty.
\end{equation}
Taking a large positive integer $N$, such that
$N>\max\{2\theta_0^{-1}, 2\beta_0^{-1}\}$.  We have
\begin{align*}
\|V_\rho\mathscr T^{2}_b(f)\|_{L^p}&\le \dsum_{d\ge 0 }\dsum_{l\in \Bbb Z } \Big(\dsum_{0\le s<Nd}+\dsum_{s>Nd}\Big)\big\|\big(\sum_{k\in \Bbb Z}
 |\big( F_{s,k,d}\Delta_{l-k}^2\big)_bf|^2\big)^{1/2}\big\|_{L^p}\\
&:=J_1+J_2.
\end{align*}
For $J_1$, by \eqref{Fsk1} and \eqref{Fsk2}, we get
\begin{align*}
J_{1}&\le
C\|b\|_{\ast}\|f\|_{L^p}\dsum_{d\ge  0}\dsum_{0\le s<Nd}\bigg(\dsum_{0\le|l|< Nd}d2^d
\sigma(E_{d})+\dsum_{|l|>Nd}2^d2^{-\theta_0|l|}\bigg) \\
&\le C\|b\|_{\ast}\|f\|_{L^p}\bigg(\dsum_{d\ge
0}d^32^d \sigma(E_{d})+\dsum_{d\ge  0 }d2^{(1-\beta_0N)d} \bigg)\\
&\le C(\|\Omega\|_{L(\log^+\!\!L)^3(\mathbf S^{n-1})}+1)\|b\|_{\ast}\|f\|_{L^p}.
\end{align*}
For $J_{2},$ using \eqref{Fsk2}, we get
\begin{align*}
J_3&\leq C\|b\|_{\ast}\dsum_{d\ge  0 }2^d\dsum_{s>Nd}2^{-\beta_0s}\bigg(\dsum_{|l|<Nd}+\dsum_{|l|\ge Nd}2^{-\theta_0|l|})\|f\|_{L^p}\\&\le C\|b\|_{\ast}\dsum_{d\ge  0 }(d2^{(1-\beta_0N)d}+2^{(1-\beta_0N-\theta_0N)d})\|f\|_{L^p}\\ &\leq C\|b\|_{\ast}\|f\|_{L^p}.
\end{align*}
Finally, combining above two estimates, we get for $1<p<\infty$
\begin{align*}\|V_\rho\mathscr T^{2}_b(f)\|_{L^p}\le
C(1+\|\Omega\|_{L(\log^+\!\!L)^3(\mathbf S^{n-1})})\|b\|_{\ast}\|f\|_{L^p}.
\end{align*}
We therefore finish the  estimate of \eqref{lnti} for $i=2$.

Now we return to prove \eqref{Fsk} and \eqref{Fsk1}.
We first prove \eqref{Fsk}. Note that $F_{s,k,d}^l=F_{s,k,d}\Delta_{l-k}$. We  write
$$
\big( F_{s,k,d}\Delta_{l-k}^2\big)_bf=
\sum_{u=0}^1 F_{s,k,d;b,u}^l\Delta_{l-k;b,1-u}f.
$$
Therefore
\begin{align}\label{fsk2}
\bigg\|\bigg(\dsum_{k\in \Bbb Z}
 |\big( F_{s,k,d}\Delta_{l-k}^2\big)_bf|^2\bigg)^{1/2}\bigg\|_{L^2}&\le
\sum_{u=0}^1\bigg\|
 \bigg(\dsum_{k\in \Bbb Z}
 |F_{s,k,d;b,u}^l\Delta_{l-k;b,1-u}f|^2\bigg)^{1/2}\bigg\|_{L^2}.
\end{align}
To proceed with the estimate \eqref{fsk2}, we define multiplier $\widetilde{F}_{s,k,d}^l$ by $\widehat{\widetilde{F}_{s,k,d}^l f}(\xi)=m_{s,k,d}^l(2^{-k}\xi)\widehat{f}(\xi)$ and $\widetilde{F}_{s,k,d;b,u}^l$ is the $u$-th commutator of $\widetilde{F}_{s,k,d}^l.$
Recall that  $m_{s,k,d}^l(\xi)=(1-\widehat{\phi_k}(\xi))\widehat{\nu_{k+s,d}}(\xi)\varphi(2^{k-l}\xi).$ Since  $supp\ (1-\widehat{\phi_k})\widehat{\nu_{k+s,d}} \subset \{\xi :|2^k\xi|>1/2\}$,
by a well-known Fourier transform estimate of Duoandikoetxea and Rubio de Francia
 (See (\cite{DR86}, p.551-552), it is easy to show that there is a  $\gamma\in(0,1)$ such that
$$|(1-\widehat{\phi_k}(\xi))\widehat{\nu_{k+s,d}}(\xi)|\le C \|\Omega_d\|_{L^\infty(\mathbf S^{n-1})} 2^{ -s\gamma}\min\{|2^k \xi|, |2^k \xi|^{-\gamma}\}$$
and
$$|\nabla[(1-\widehat{\phi_k}(\xi))\widehat{\nu_{k+s,d}}(\xi)]|\le C\|\Omega_d\|_{L^\infty(\mathbf S^{n-1})} 2^k 2^s.$$
As a result, we have the following estimates\begin{equation}\label{mskd0}supp\, m_{s,k,d}^l(2^{-k}\xi)\subset\{|\xi|\le 2^l\}\end{equation}
\begin{equation}\label{mskdl}
|m_{s,k,d}^l(2^{-k}\xi)|\le C\|\Omega_d\|_{L^\infty(\mathbf S^{n-1})} 2^{ -\gamma s}  \min\{2^l, 2^{-\gamma l}\},
\end{equation}
and
\begin{equation}\label{mskd2}|\nabla m_{s,k,d}^l(2^{-k}\xi)|\le C\|\Omega_d\|_{L^\infty(\mathbf S^{n-1})}  2^s.\end{equation}
Applying Lemma \ref{mulcom} with $\sigma=2^l$ to \eqref{mskd0}-\eqref{mskd2}, there exist  constants $\beta\in (0,1)$ and $\theta\in (0,1)$ such that
$$ \|\widetilde{F}_{s,k,d;b,u}^l f\|_{L^2} \le C\|\Omega_d\|_{L^\infty(\mathbf S^{n-1})}\|b\|_{\ast}^u2^{-\beta s}2^{-\theta|l|}\|f\|_{L^2},\ \ \ \ \ for \ \ \ l\in \Bbb Z \ \ and\ \ s\ge 0.$$
Further, the dilation-invariance implies
\begin{equation}\label{Fskdb2}
\|F_{s,k,d;b,u}^lf\|_{L^2}\le C\|\Omega_d\|_{L^\infty(\mathbf S^{n-1})}\|b\|_{\ast}^u 2^{-\beta s}2^{-\theta|l|}\|f\|_{L^2},\ \ \ \ \ \text{for} \ l\in \Bbb Z\  \text{and}\ s\ge 0.
\end{equation}
Then by \eqref{Fskdb2} and Lemma \ref{Deltacom}, we get
\begin{align*}
\bigg\|\bigg(\dsum_{k\in \Bbb Z}
 |\big( F_{s,k,d}\Delta_{l-k}^2\big)_bf|^2\bigg)^{1/2}\bigg\|_{L^2}&\le C\sum_{u=0}^1\|\Omega_d\|_{L^\infty(\mathbf S^{n-1})}2^{-\beta s}2^{-\theta |l|}\|b\|_{\ast}^{u} \bigg\|
 \bigg(\dsum_{k\in \Bbb Z}
 |\Delta_{l-k;b,1-u}f|^2\bigg)^{1/2}\bigg\|_{L^2}\\
 &\le C\|\Omega_d\|_{L^\infty(\mathbf S^{n-1})}2^{-\beta s}2^{-\theta |l|}\|b\|_{\ast}\|f\|_{L^2}.
\end{align*}
This gives \eqref{Fsk}.
Secondly,  we turn to prove \eqref{Fsk1}. Write
\begin{align*}
\big( F_{s,k,d}\Delta_{l-k}^2\big)_bf(x)=\sum_{u=0}^1
 F_{s,k,d;b,u} \Delta_{l-k;b,1-u}^2f.
\end{align*}
 By the Minkowski inequality, we get
\begin{align}\label{gsp}\bigg\|\bigg(\dsum_{k\in \Bbb Z}
 |\big( F_{s,k,d}\Delta_{l-k}^2\big)_bf|^2\bigg)^{1/2}\bigg\|_{L^p}&\le\sum_{u=0}^1\bigg\|
 \bigg(\dsum_{k\in \Bbb Z}
 |F_{s,k,d;b,u}\Delta_{l-k;b,1-u}^2f|^2\bigg)^{1/2}\bigg\|_{L^p}.
\end{align}
To proceed with the above estimate, we need the following lemma,
 which can be proved as that in (\cite[p.\,544]{DR86}).

\begin{lemma}\label{scom} Suppose that $\{\sigma_j\}_{j\in\mathbb{Z}}$ is a
sequence of finite Borel measures, $T_jf=\sigma_j\ast f$ and ${\overline{T}_j}f=|\sigma_j|\ast f$ for any $j\in \Bbb Z$. Let $b\in BMO$ and $T_{j;b,u}f(x)=T_j((b(x)-b(\cdot))^uf)(x),\,\,u\in \Bbb N.$ If the maximal operator
$T_{b,u}^{*}(f)=\displaystyle\sup_{j\in \Bbb Z}|{\overline{T}}_{j;b,u}f|$ is bounded on
$L^{p_0}$ for any fixed $1<p_0<\infty$, then for any $1<p<\infty,$
$$\bigg\|\bigg(\dsum_{j\in \Bbb Z}|T_{j;b,u}g_j|^2\bigg)^{1/2}\bigg\|_{L^p}\le C\|b\|_{\ast}^u\bigg\|\bigg(\dsum_{j\in \Bbb Z}|g_j|^2\bigg)^{1/2}\bigg\|_{L^p}. $$
\end{lemma}
Recall that $
 F_{s,k,d}f(x)=(\delta_0-\phi_k)\ast\nu_{k+s}\ast f(x).$ Define by $
 \overline{F}_{s,k,d}f(x)=|(\delta_0-\phi_k)\ast\nu_{k+s}|\ast f(x).$
Since  for some $1<p_0<\infty$,
\begin{align*}
 \|\sup_{k\in \Bbb Z} |\overline{F}_{s,k,d;b,1}f|\|_{L^{p_0}}\le
 \|{\Omega}_d\|_{L{\log^+L}(\mathbf S^{n-1})}\|b\|_{\ast}\|
f\|_{L^{p_0}}
\end{align*}
and \begin{align*}
 \|\sup_{k\in \Bbb Z} |\overline{F}_{s,k,d}f|\|_{L^{p_0}}\le
 \|{\Omega_d}\|_{L^1(\mathbf S^{n-1})}\|
f\|_{L^{p_0}}
\end{align*}(see \cite{H, G}). So by \eqref{gsp},  Lemma \ref{scom}
 and Lemma \ref{Deltacom}, we get for  $1<p<\infty$,
\begin{align}\nonumber\bigg\|\bigg(\dsum_{k\in \Bbb Z}
 |\big( F_{s,k,d}\Delta_{l-k}^2\big)_bf|^2\bigg)^{1/2}\bigg\|_{L^p}\le
 C\|\Omega_d\|_{L{\log^+L}(\mathbf S^{n-1})}\|b\|_{\ast}\|f\|_{L^p},
\end{align} which gives \eqref{Fsk1}.


{{ \textbf{Estimate of \eqref{lnti} for $i=3$}}.}
 We have the following pointwise estimate
\begin{align}
\nonumber V_\rho\mathscr T_b^3(f)(x)
&\le\dsum_{s<0}\big(\dsum_{k\in \Bbb Z}\big|(\phi_k\ast\nu_{k+s})_b
f(x)\big|^2\big)^{1/2}.
\end{align}

The proofs are essentially similar to  the proof of  (\ref{lnti}) for $i=2$. More precisely,  we need to give the estimates on  the left hand side of (\ref{Fsk})-(\ref{Fsk1}) with replacing $(\delta_0-\phi_k)\ast\nu_{k+s}$ by $\phi_k\ast\nu_{k+s}$. Since  $supp\ \widehat{\phi_k}\widehat{\nu_{k+s}} \subset \{\xi :|2^k\xi|<1\}$ and $\Omega$ satisfies \eqref{can of O}, then it is easy to see that
$$|\widehat{\phi_k\nu_{k+s}} (\xi)|\le C 2^{s}\|\Omega\|_{L^1(\mathbf S^{n-1})} \min\{|2^k \xi|, |2^k \xi|^{-1}\}$$  and
\begin{align*}
|\nabla\widehat{\nu_{k+s}} (\xi)|&\le C2^{(k+s)}\|\Omega\|_{L^1(\mathbf S^{n-1})}.
\end{align*}
 Set   $$R_{s,k}(\xi)=\widehat{\phi_k}(\xi)\widehat{\nu_{k+s}}(\xi),\,\,\,\,\,R_{s,k}^l(\xi)=R_{s,k}(\xi)\varphi(2^{k-l}\xi).$$
Using the two above inequalities, we have the following estimate \begin{align*}\label{m2}
supp R_{s,k}^l(2^{-k}\xi)\subset\{|\xi| \le 2^l\},
\end{align*}
\begin{align*}
|R_{s,k}^l(2^{-k}\xi)|\le C 2^{ s} \min\{2^{l}, 2^{- l}\}\|\Omega\|_{L^1(\mathbf S^{n-1})}
\end{align*} and
 \begin{align*}|\nabla (R_{s,k}^l(2^{-k}\xi))|&\le  C2^s\|\Omega\|_{L^1(\mathbf S^{n-1})}.\end{align*}
Then apply Lemma \ref{mulcom} and the same arguments of the proofs of  (\ref{lnti}) for $i=2$, then  the right hand side of (\ref{Fsk}) is controlled by $ C2^{  s}2^{-\theta |l|}\|\Omega\|_{L^1(\mathbf S^{n-1})}\| b\|_{\ast}\|f\|_{L^2}$ for $\theta>0.$
It is also  easy to get the same estimates in the right hand side of (\ref{Fsk1}) by using Lemma \ref{Deltacom}, Lemma  \ref{scom} and the results in \cite{H, G}.
Then we get for $1<p<\infty$
$$\| V_\rho\mathscr T_b^3(f)\|_{L^p}\le
C\|\Omega\|_{L{\log^+L}(\mathbf S^{n-1})}\| b\|_{\ast}\|f\|_{L^p}.
$$
We therefore finish the  proof of \eqref{lnti} for $i=3$.\qed

\subsection{Proof of \eqref{short variation singular}}

In the section, we give the proof of \eqref{short variation singular}.
For $t\in[1,2)$, we define $\nu_{0,t}$ as
$$\nu_{0,t}(x)=\frac{\Omega(x')}{|x|^n}\chi_{\{t\leq|x|\leq2\}}(x)$$
and $\nu_{j,t}(x)={2^{-jn}}\nu_{0,t}(2^{-j}x)$ for $j\in\mathbb{Z}$. Denote $T_{j,t}$ by $T_{j,t}f(x)=\nu_{j,t}\ast f(x)$ and $T_{j,t;b}$ by $T_{j,t;b}f(x)=T_{j,t}(b(x)-b(\cdot))f)(x)$.
 Observe that
$V_{2,j}(\mathcal{T}_bf)(x)$ is just the strong  $2$-variation function of the family
$\{T_{j,t;b} f(x)\}_{t\in[1,2)}$ and $\sum_{k\in\Bbb Z}\Delta_k^2=\mathcal{I}$\,(identity operator), hence
\begin{align*}
S_{2}(\mathcal{T}_bf)(x)&=\Big(\sum_{j\in\mathbb{Z}}|V_{2,j}(\mathcal{T}_bf)(x)|^2\Big)^{\frac{1}{2}}
=\Big(\sum_{j\in\mathbb{Z}}\|\{T_{j,t;b} f(x)\}_{t\in[1,2)}\|_{V_2}^2\Big)^{\frac{1}{2}}\\
&\leq\sum_{k\in\mathbb{Z}}\Big(\sum_{j\in\mathbb{Z}}\|\{\big(T_{j,t}\Delta_{k-j}^2\big)_b
f(x)\}_{t\in[1,2)}\|_{V_2}^2\Big)^{\frac{1}{2}}\\&:=
\sum_{k\in\mathbb{Z}}S_{2,k}(\mathcal{T}_bf)(x).
\end{align*}
For $d\ge 0$, set
   $$\nu_{j,t,d}(x)=\frac{\Omega_d(x)}{|x|^{n}}\chi_{\{2^jt\le |x|< 2^{j+1}\}}(x).$$
$T_{j,t,d}$ is defined as $T_{j,t}$ by replacing $\nu_{j,t}$ by $\nu_{j,t,d}$. Decompose $\Omega$ as  in the estimate of \eqref{lnti} for $i=2$. Then,
$$
S_{2}(\mathcal{T}_bf)(x)\le\dsum_{d\ge 0 }\dsum_{k\in \Bbb Z}\Big(\sum_{j\in\mathbb{Z}}\|\big\{\big(T_{j,t,d}\Delta_{k-j}^2\big)_b
f(x)\big\}_{t\in[1,2)}\|_{V_2}^2\Big)^{\frac{1}{2}}:=\dsum_{d\ge 0 }\dsum_{k\in \Bbb Z}
S_{2,k,d}(\mathcal{T}_bf)(x).
$$
\begin{proposition}\label{S2k} For $d\ge0$ and $k\in\mathbb Z$,
following conclusions hold:\\
{\rm (i)}\quad  There exist  a constant $
\theta\in(0,1)$ such that
\begin{equation}\label{S2k2}
\|S_{2,k,d}(\mathcal{T}_bf) \|_{L^2}\le
C\|\Omega_d\|_{L^\infty(\mathbf S^{n-1})}2^{-\frac{\theta}{2}|k|}\|
b\|_{\ast}\|f\|_{L^2};
\end{equation}
{\rm (ii)}\quad For $1<p<\infty,$
\begin{equation}\label{S2kp}
\|S_{2,k,d}(\mathcal{T}_bf) \|_{L^p}\le C\|\Omega_d\|_{L\log^{+}L(\mathbf S^{n-1})} \|
b\|_{\ast}\|f\|_{{L}^p}.
\end{equation}
The constants $C's$ in \eqref{S2k2} and  \eqref{S2kp} are
independent of $k$.
\end{proposition}
  Using the same argument of \eqref{lnti} for $i=2,$ we may
 finish the proof of \eqref{short variation singular} by using Proposition \ref{S2k}. We omit the details.

{\emph{Proof of Proposition  \ref{S2k}.}} Without loss of generality, we will use $S_{2,k}$ to replace with $S_{2,k,d}.$
To deal with \eqref{S2k2}, we borrow the fact $\|\mathfrak{a}\|_{V_2}\le\|\mathfrak{a}\|_{L^2}^{1/2}\|\mathfrak{a}'\|_{L^2}^{1/2}$, where $\mathfrak{a}'=\{\frac{d}{dt}a_t:t\in\mathbb R\}$.
It is a special case of (39) in \cite{JSW08}.  Then,
\begin{align*}
[S_{2,k}(\mathcal{T}_bf)(x)]^2&\leq\sum_{j\in\mathbb{Z}}\bigg(\dint_{1}^{2}|\big(T_{j,t}
\Delta_{k-j}^2\big)_bf(x)|^2\frac{dt}{t}\bigg)^{1/2}\bigg(\dint_{1}^{2}|\big(\frac{d}{dt}T_{j,t}
\Delta_{k-j}^2\big)_bf(x)|^2\frac{dt}{t}\bigg)^{1/2}\\
&:=\sum_{j\in\mathbb{Z}}I_{1,k}f(x)\cdot I_{2,k}f(x).
\end{align*}
By the Cauchy-Schwarz inequality, we have
\begin{align*}\nonumber
&\big\|S_{2,k}(\mathcal{T}_bf)\big\|^2_{L^2}\le
\|I_{1,k}f\|_{L^2}\|I_{2,k}f\|_{L^2}.
\end{align*}
We estimate $\|I_{1,k}f\|_{L^2}$ and $\|I_{2,k}f\|_{L^2}$, respectively.
To estimate $\|I_{1,k}f\|_{L^2}$, we need the following
estimates: for some  $\gamma>0,$
\begin{equation*}
|\widehat{\nu_{j,t}}(\xi)|\leq C\|\Omega\|_{L^\infty(\mathbf S^{n-1})}\min\{|2^j\xi|^{-\gamma},
       2^j|\xi|\}
\end{equation*}
and
\begin{equation*}
|\nabla\widehat{\nu_{j,t}}(\xi)|\leq C\|\Omega\|_{L^\infty(\mathbf S^{n-1})}
\end{equation*}
uniformly in $t\in[1,2)$, which have been essentially proved in \cite{DR86} and \cite{H}.
Similarly to the proof of \eqref{Fsk}, we get for some $v>0,$
\begin{equation}\label{I1kd2}
\|I_{1,k}f\|_{L^2}\le C\|\Omega\|_{L^\infty(\mathbf S^{n-1})}\|b\|_{\ast}\min\{ 2^{ vk}, 2^{-\gamma v k}\}\|f\|_{L^2}.
\end{equation}
Next, we estimate $\|I_{2,k}f\|_{L^2}$.  Write
\begin{align*}
\big(\frac{d}{dt}T_{j,t}
\Delta_{k-j}^2\big)_bf&=\Delta_{k-j;b}\dfrac{d}{dt}T_{j,t}\Delta_{k-j} f+
\Delta_{k-j} \big(\dfrac{d}{dt}T_{j,t}\big)_b\Delta_{k-j}f+\Delta_{k-j}\dfrac{d}{dt}T_{j,t}  \Delta_{k-j;b}f.
\end{align*}
Thus we get
\begin{align*}
\|I_{2,k}f\|_{L^2}&\le \bigg(\dint_{1}^{2}\bigg\|\bigg(\dsum_{j\in \Bbb Z}|\Delta_{k-j;b}\dfrac{d}{dt}T_{j,t}
\Delta_{k-j}f|^2\bigg)^{1/2}\bigg\|_{L^2}^2\frac{dt}{t}\bigg)^{1/2}\\&+\bigg(\dint_{1}^{2}\bigg\|\bigg(\dsum_{j\in \Bbb Z}| \Delta_{k-j}\big(\dfrac{d}{dt}T_{j,t}\big)_b
\Delta_{k-j}f|^2\bigg)^{1/2}\bigg\|_{L^2}^2\frac{dt}{t}\bigg)^{1/2}\\&+\bigg(\dint_{1}^{2}\bigg\|\bigg(\dsum_{j\in \Bbb Z}| \Delta_{k-j}\dfrac{d}{dt}T_{j,t}
\Delta_{k-j;b}f|^2\bigg)^{1/2}\bigg\|_{L^2}^2\frac{dt}{t}\bigg)^{1/2}\\&:=I+II+III.
\end{align*}
To estimate $I,\,II,\,III$, respectively,
we need the following elementary fact
\begin{align}
\nonumber\frac{d}{dt}T_{j,t}h(x)&=\dfrac{d}{dt}\bigg[\int_{2^jt<|y|\le 2^{j+1}}\dfrac{\Omega(y')}{|y|^n}h(x-y)dy\bigg]\\
\nonumber&=\dfrac{d}{dt}\bigg[\dint_{\mathbf S^{n-1}}\Omega(y')\dint_{2^jt}^{2^{j+1}}\frac{1}{r}h(x-ry')drd\sigma(y')\bigg]\\
\label{dTjtd}&=-\dfrac{1}{t}\dint_{\mathbf S^{n-1}}\Omega(y')h(x-2^jty')d\sigma(y')
\end{align}
 and \begin{equation}\label{Tjtdsq2}
\|T_{j,t}^*h\|_{L^2}\le C\|\Omega\|_{L^1(\mathbf S^{n-1})}\|h\|_{L^2},
\end{equation} where $$T_{j,t}^*h(x)=\int_{\mathbf S^{n-1}}|\Omega(y')||h(x-2^jty')|d\sigma(y')$$ for  $t\in [1,2).$
We now estimate $I$.
Indeed, by  \eqref{Tjtdsq2}, Lemma \ref{Deltacom} and  Littlewood-Paley theory, we
have
\begin{align*}
I
&\le C\|b\|_{\ast}\bigg(\dint_{1}^{2} \dsum_{j\in
\Bbb Z}\big\|T_{j,t}^*
\Delta_{k-j}f
\big\|_{L^2}^2\frac{dt}{t}\bigg)^{1/2}\\
&\le C\|b\|_{\ast}\|\Omega\|_{L^1(\mathbf S^{n-1})}\bigg\|\bigg(\dsum_{j\in
\Bbb Z}\big|
\Delta_{k-j}f\big|^2
\bigg)^{1/2}\bigg\|_{L^2}\\
&\le C\|b\|_{\ast}\|\Omega\|_{L^1(\mathbf S^{n-1})}\|f\|_{L^2}.
\end{align*}
Similarly, by  \eqref{Tjtdsq2}, Lemma \ref{Deltacom} and  Littlewood-Paley theory, we get
$$III\le C\|b\|_{\ast}\|\Omega\|_{L^1(\mathbf S^{n-1})}\|f\|_{L^2}.$$
For $II,$ we will apply  the Bony paraproduct to do this. Let $\varpi\in {\mathscr S}({\Bbb R}^n)$ be a
radial function satisfying $0\le \varpi\le1$ with its support is in
the unit ball and $\varpi(\xi)=1$ for $|\xi|\le \frac12$. The
function $\psi(\xi)=\varpi(\frac\xi2)-\varpi(\xi)\in {\mathscr
S}({\Bbb R}^n)$ supported on $\{\frac12\le|\xi|\le 2\}$ and
satisfies the identity $\sum_{j\in\mathbb
Z}\psi(2^{-j}\xi)=1$ for $\xi\neq0$.
For $j\in\mathbb Z$, denote by $\Theta_j$ and $G_j$ the convolution
operators whose the symbols are $\psi(2^{-j}\xi)$ and
$\varpi(2^{-j}\xi)$, respectively. That is, $\Theta_j$ and $G_j$
are defined by $\widehat{\Theta_jf}(\xi)
=\psi(2^{-j}\xi)\hat{f}(\xi)$ and
$\widehat{G_jf}(\xi)=\varpi(2^{-j}\xi)\hat{f}(\xi)$
(see [\ref{G}]). The
paraproduct of Bony [\ref{B}] between two functions $f$, $g$ is
defined by
$$\pi_f(g)=\dsum_{j\in\mathbb Z}(\Theta_jf)(G_{j-3}g).$$
At least formally, we have the following Bony decomposition
\begin{equation}\label{Bd}
fg=\pi_f(g)+\pi_g(f)+R(f,g)\quad\text{with}\quad
R(f,g)=\dsum_{i\in\mathbb Z}\dsum_{|k-i|\le 2}(\Theta_if)(\Theta_k
g).
\end{equation}
Denote $f_{k,j}:=\Delta_{k-j}f$ and $T_{j,t}'f:=\frac{d}{dt}T_{j,t}f$. By \eqref{Bd},  we have
\begin{align*}
[b,T_{j,t}']f_{k,j}(x)
&=b(x)(T_{j,t}'f_{k,j})(x)-T_{j,t}'(bf_{k,j})(x)\\
&=[\pi_{T_{j,t}'f_{k,j}}(b)(x)-T_{j,t}'\big(\pi_{f_{k,j}}(b)\big)(x)]\\
&+[R(b,T_{j,t}'f_{k,j})(x)-T_{j,t}'\big(R(b,f_{k,j})\big)(x)]\\
&+ [\pi_b(T_{j,t}'f_{k,j})(x)-T_{j,t}'
\big(\pi_b(f_{k,j})\big)(x)].
\end{align*}
 Thus
\begin{align*}II&\le
\bigg(\dint_{1}^{2}\bigg\|\bigg(\sum_{j\in\mathbb{Z}}|\Delta_{k-j}[\pi_{T_{j,t}'f_{k,j}}(b)-T_{j,t}'\big(\pi_{f_{k,j}}(b)\big)]|^2\bigg)^{1/2}\bigg\|_{L^2}^2\frac{dt}{t}\bigg)^{1/2}
\\
&+\bigg(\dint_{1}^{2}\bigg\|\bigg(\sum_{j\in\mathbb{Z}}|\Delta_{k-j}[R(b,T_{j,t}'f_{k,j})-T_{j,t}'\big(R(b,f_{k,j})\big)]|^2\bigg)^{1/2}\bigg\|_{L^2}^2\frac{dt}{t}\bigg)^{1/2}
\\
&\quad+\bigg(\dint_{1}^{2}\bigg\|\bigg(\sum_{j\in\mathbb{Z}}|\Delta_{k-j} [\pi_b\big(T_{j,t}'f_{k,j}\big)-T_{j,t}'
\big(\pi_b(f_{k,j})\big)]|^2\bigg)^{1/2}\bigg\|_{L^2}^2\frac{dt}{t}\bigg)^{1/2}\\
&:=II_1+II_2+II_3.\end{align*}
We will estimate $II_i,\,i=1,2,3,$ respectively.
For $II_1$,  note that $\Theta_i\Delta_{k-j}g=0$ for $g\in
{\mathscr S}'({\Bbb R}^n)$
 when $|i-(k-j)|\ge 3$, by \eqref{dTjtd}, we
get
\begin{align}
\nonumber&\big[\pi_{T_{j,t}'f_{k,j}}(b)(x)-T_{j,t}'\big(\pi_{f_{k,j}}(b)\big)(x)\big]\\
\nonumber=&\dsum_{i\in \Bbb Z}\bigg\{(T_{j,t}'\Theta_{i}\Delta_{k-j}f)(x)(G_{i-3}b)(x)-
T_{j,t}'[(\Theta_i\Delta_{k-j}f)(G_{i-3}b)](x)\bigg\}\\
\label{Tjtdcom}=&\dsum_{|i-(k-j)|\le 2}\dfrac{1}{t}\dint_{\mathbf S^{n-1}}\Omega(y')(G_{i-3}b(x)-G_{i-3}b(x-2^jty'))(\Theta_i\Delta_{k-j}f)(x-2^jty')d\sigma(y').
\end{align}

To estimate the above inequality, we need the follow lemma.

\begin{lemma}\label{difGkb} (\cite{CD1}) For any fixed $0<\tau<1/2$, we have
$$|G_kb(x)-G_kb(y)|\le
C\frac{2^{k\tau}}{\tau}|x-y|^\tau\|b\|_{\ast},$$ where $C$ is independent of $k$ and
$\tau$.
\end{lemma}
By Lemma \ref{difGkb},  note that for $t\in [1,2)$,  we get
\begin{align}
\nonumber&\Big|\dfrac{1}{t}\dint_{\mathbf S^{n-1}}\Omega(y')(G_{i-3}b(x)-G_{i-3}b(x-2^jty'))(\Theta_i\Delta_{k-j}f)(x-2^jty')d\sigma(y')\Big|\\&\nonumber
\le C\dfrac{2^{(i+j)\tau}}{\tau}\|b\|_{\ast}\dint_{\mathbf S^{n-1}}|\Omega(y')||(\Theta_i\Delta_{k-j}f)(x-2^jty')|d\sigma(y')\\
\label{cGiTjtd}&\le C\dfrac{2^{(i+j)\tau}}{\tau}\|b\|_{\ast}T_{j,t}^*(\Theta_i\Delta_{k-j}f)(x).
\end{align}
Then by  \eqref{Tjtdcom}-\eqref{cGiTjtd},
\eqref{Tjtdsq2} and   Littlewood-Paley theory,
 we have
\begin{align*}
II_1&\le C\|b\|_{\ast}\dfrac{2^{k\tau}}{\tau}\dsum_{|l|\le
2}\bigg(\int_1^2\Big\|\bigg(\dsum_{j\in\mathbb
Z}|T_{j,t}^*(\Theta_{k-j+l}\Delta_{k-j}f)|^2\bigg)^{1/2}
\Big\|_{L^2}^2\dfrac{dt}{t}\bigg)^{1/2}\\
&\le C\dfrac{2^{k\tau}}{\tau}\|b\|_{\ast}\|\Omega\|_{L^1(\mathbf S^{n-1})}\|f\|_{ L^2},
\end{align*}
where $C$ is
independent of $k$ and $\tau$.

Next, we estimate $II_2$. Clearly,
$\Theta_{i+l}\Delta_{k-j}g=0$ for $g\in \mathscr{S}^{'}({\Bbb R}^n)$ when
$|l|\le2$ and $|i-(k-j)|\ge 8$. Thus by \eqref{dTjtd},
\begin{align}
\nonumber&[R(b,T_{j,t}'f_{k,j})(x)-T_{j,t}'\big(R(b,f_{k,j})\big)(x)]\\
\nonumber=&\dsum_{i\in\mathbb Z}\dsum_{|l|\le 2}
(\Theta_ib)(x)(T_{j,t}'\Theta_{i+l}\Delta_{k-j}f)(x)-
T_{j,t}'\Big(\dsum_{i\in\mathbb Z}
\dsum_{|l|\le 2}(\Theta_ib)(\Theta_{i+l}\Delta_{k-j}f)\Big)(x)\\
\nonumber=&\dsum_{l=-2}^2\dsum_{|i-(k-j)|\le 7}\bigg((\Theta_ib)(x) (T_{j,t}'
\Theta_{i+l}\Delta_{k-j}f)(x)-
T_{j,t}'\big((\Theta_ib)(\Theta_{i+l}\Delta_{k-j}f)\big)(x)\bigg)\\
\label{RbTjtd}=&\dsum_{l=-2}^2\dsum_{|i-(k-j)|\le 7}\dfrac{1}{t}\dint_{\mathbf S^{n-1}}\Omega(y')(\Theta_ib(x)-\Theta_ib(x-2^jty'))\Theta_{i+l}\Delta_{k-j}f(x-2^jty')d\sigma(y').
\end{align}
Note that for $t\in [1,2)$, by $\sup\limits_{i\in\mathbb
Z}\|\Theta_i(b)\|_{L^\infty}\le C\|b\|_{\ast}$ (see \cite{G}), we get for any $h\in L^p(\Bbb R^n),$
\begin{align}
\nonumber&\Big|\dfrac{1}{t}\dint_{\mathbf S^{n-1}}\Omega(y')(\Theta_ib(x)-\Theta_ib(x-2^jty'))h(x-2^jty')d\sigma(y')\Big|\\
\nonumber&\le C\sup\limits_{i\in\mathbb
Z}\|\Theta_i(b)\|_{L^\infty}\dint_{\mathbf S^{n-1}}|\Omega(y')||h(x-2^jty')|d\sigma(y')\\&\le\label{TibTjtd}
C\|b\|_{\ast}T_{j,t}^*h(x).
\end{align}
Using \eqref{RbTjtd}-\eqref{TibTjtd},
\eqref{Tjtdsq2} and  Littlewood-Paley theory, we have
\begin{align*}
II_2&\le C\|b\|_{\ast}\|\Omega\|_{L^1(\mathbf S^{n-1})}\dsum_{|l|\le
7}\bigg(\dint_{1}^{2}\big\|\big(\sum_{j\in\mathbb{Z}}|T_{j,t}^*(\Theta_{k-j+l}\Delta_{k-j}f)|^2\big)^{1/2}\big\|_{L^2}^2\frac{dt}{t}\bigg)^{1/2}\\
&\le C\|b\|_{\ast}\|\Omega\|_{L^1(\mathbf S^{n-1})}\|f\|_{ L^2}.
\end{align*}

Finally, we estimate $II_3$.  Note that $\Delta_{k-j}\big((\Theta_i g)(
G_{i-3}h)\big)=0$ for $g,\, h\in \mathscr{S}^{'}({\Bbb R}^n)$ if
$|i-(k-j)|\ge 5$. Thus we get
\begin{align*}
&\Delta_{k-j}[\pi_b\big(T_{j,t}'f_{k,j}\big)(x)-T_{j,t}'
\big(\pi_b(f_{k,j})\big)](x)\\
=&\Delta_{k-j}\Big(\dsum_{i\in\mathbb
Z}(\Theta_ib)(G_{i-3}T_{j,t}'\Delta_{k-j}f)-T_{j,t}'\big(\dsum_{i\in\mathbb
Z}(\Theta_ib)
(G_{i-3}\Delta_{k-j}f)\big)\Big)(x)\\
=&\dsum_{|i-(k-j)|\le 4}\Delta_{k-j}\Big\{\big(\Theta_ib T_{j,t}'\big(G_{i-3}
\Delta_{k-j}f)\big)-T_{j,t}'
\big((\Theta_ib)(G_{i-3}\Delta_{k-j}f)\big)\Big\}(x).
\end{align*}
Since \begin{align*}
&\big(\Theta_ib T_{j,t}'\big(G_{i-3}
\Delta_{k-j}f)\big)-T_{j,t}'
\big((\Theta_ib)(G_{i-3}\Delta_{k-j}f)\big)(x)\\&=\dfrac{1}{t}\dint_{\mathbf S^{n-1}}\Omega(y')(\Theta_ib(x)-\Theta_ib(x-2^jty'))G_{i-3}\Delta_{k-j}f(x-2^jty')d\sigma(y').
\end{align*}
Thus, by \eqref{TibTjtd}, \eqref{Tjtdsq2} and  Littlewood-Paley theory,
 we get
\begin{align*}
II_3&\le C\|b\|_{\ast}\|\Omega\|_{L^1(\mathbf S^{n-1})}\dsum_{|l|\le
4}\bigg(\dint_{1}^{2}\bigg\|\bigg(\sum_{j\in\mathbb{Z}}|T_{j,t}^*G_{k-j+l-3}\Delta_{k-j}f|^2\bigg)^{1/2}\bigg\|_{L^2}^2\frac{dt}{t}\bigg)^{1/2}\\
&\le
C\|b\|_{\ast}\|\Omega\|_{L^1(\mathbf S^{n-1})}\|f\|_{{L}^2}.
\end{align*}
Together with the estimates of $II_1,$ $II_2$ and $II_3,$ we  get
\begin{equation}\label{L2t}
II\le C\max\{2,\frac{2^{\tau
k}}{\tau}\}\|b\|_{\ast}\|\Omega\|_{L^1(\mathbf S^{n-1})}\|f\|_{{L}^2}\quad\hbox{ for}\quad
k\in \Bbb Z,
\end{equation}
where $C$ is independent of $k$ and $\tau$.
Taking $\tau=\frac{1}{|k|}$ in \eqref{L2t}, we get
$$II\le C(|k|+1)\|b\|_{\ast}\|\Omega\|_{L^1(\mathbf S^{n-1})}\|f\|_{{L}^2}\quad\hbox{ for}\quad
k\in \Bbb Z.$$
Combining this with the estimate of $I$ and $III$, we get
\begin{equation}\label{I2kd}
\|I_{2,k}f\|_{L^2}\le  C(|k|+1)\|b\|_{\ast}\|\Omega\|_{L^1(\mathbf S^{n-1})}\|f\|_{{L}^2}.
\end{equation}
Combining the estimates of \eqref{I1kd2} and  \eqref{I2kd}, we get for some constant $\theta\in (0,1)$ and  $k\in \Bbb Z,$
\begin{align*}
\big\|S_{2,k}(\mathcal{T}_bf)\big\|^2_{L^2}&\le
C\min\{ 2^{ vk}, 2^{-\gamma v k}\}(1+|k|)\|b\|_{\ast}^2\|\Omega\|_{L^\infty(\mathbf S^{n-1})}^2\|f\|_{ L^2}^2\\&\le
C2^{-\theta|k|}\|b\|_{\ast}^2\|\Omega\|_{L^\infty(\mathbf S^{n-1})}^2\|f\|_{ L^2}^2.
\end{align*} This finishes the proof of \eqref{S2k2}.

\emph{Proof of  \eqref{S2kp}}. Let $$
B=\Big\{\{a_{j,t}\}_{j\in \Bbb Z,\, t\in [1,2)}:\|{a_{j,t}}\|_B:=\big(\sum_{j\in \Bbb Z}\|{a_{j,t}}\|^2_{V_2}\big)^{1/2}<\infty\Big\}.
$$
Clearly, $(B,\|\cdot\|_B)$ is a Banach space.Then,
\begin{align*}
S_{2,k}(\mathcal{T}_bf)(x)&=\Big(\dsum_{j\in\mathbb{Z}}\dsup_{\substack
{t_1<\cdots<t_N\\
[t_l,t_{l+1}]\subset[1,2)}}\dsum_{l=1}^{N-1}\big|\big(T_{j,t_l}\Delta_{k-j}^2\big)_b
f(x)-\big(T_{j,t_{l+1}}\Delta_{k-j}^2\big)_b
f(x)|^2\Big)^{\frac{1}{2}}\\
&=\Big(\dsum_{j\in\mathbb{Z}}\dsup_{\substack
{t_1<\cdots<t_N\\
[t_l,t_{l+1}]\subset[1,2)}}\dsum_{l=1}^{N-1}\big|\big(T_{j,t_l,t_{l+1}}\Delta_{k-j}^2\big)_b
f(x)|^2\Big)^{\frac{1}{2}},
\end{align*} where $$T_{j,t_l,t_{l+1}}f(x)=\dint_{2^jt_l<|y|\le 2^jt_{l+1}}f(x-y)\frac{\Omega(y)}{|y|^n}dy\,\,\, \hbox{and}\,\,\,\,[t_l,t_{l+1}]\subset[1,2).$$
Then we get
\begin{align*}
S_{2,k}(\mathcal{T}_bf)(x)&\le  \sum_{u=0}^1\Big(\dsum_{j\in\mathbb{Z}}\dsup_{\substack
{t_1<\cdots<t_N\\
[t_l,t_{l+1}]\subset[1,2)}}\dsum_{l=1}^{N-1}\big|T_{j,t_l,t_{l+1};b,u}\Delta_{k-j;b,1-u}^2
f(x)|^2\Big)^{\frac{1}{2}}.
\end{align*}
For $[t_l,t_{l+1}]\subset[1,2)$ and $u\in\Bbb N$, let  $${\widetilde{T}}_{j,t_l,t_{l+1};  b,u}f(x)=\dint_{2^jt_l<|x-y|\le 2^jt_{l+1}}|f(y)||b(x)-b(y)|^u\frac{|\Omega(x-y)|}{|x-y|^n}dy.$$
Then,
\begin{align*}
S_{2,k}(\mathcal{T}_bf)(x)&\le \sum_{u=0}^1\Big(\dsum_{j\in\mathbb{Z}}\dsup_{\substack
{t_1<\cdots<t_N\\
[t_l,t_{l+1}]\subset[1,2)}}\dsum_{l=1}^{N-1}\big|{\widetilde{T}}_{j,t_l,t_{l+1};  b,u}\Delta_{k-j;b,1-u}^2
f(x)|^2\Big)^{\frac{1}{2}}\\&= \sum_{u=0}^1\Big(\dsum_{j\in\mathbb{Z}}\dsup_{\substack
{t_1<t_N\\
[t_1,t_{N}]\subset[1,2)}}\big|\widetilde{T}_{j,t_1,t_{N};b,u}(\Delta_{k-j;b,1-u}^2
f)(x)|^2\Big)^{\frac{1}{2}}.
\end{align*}
 Therefore, we get
 \begin{equation}\label{S2kdm}
 S_{2,k}(\mathcal{T}_bf)(x)
 \le C\sum_{u=0}^1\Big(\dsum_{j\in\mathbb{Z}}\big|T^*_{\Omega;b,u}(\Delta_{k-j;b,1-u}^2
f)(x)|^2\Big)^{\frac{1}{2}},\end{equation} where $$T^*_{\Omega;b,u}f(x)=\sup_{r>0}\frac{1}{r^n}\int_{|x-y|<r}|f(y)||b(x)-b(y)|^u|\Omega(x-y)|\,dy.$$
Since for $1<p<\infty$ and $u\in\{0,1\}$,
$$
\bigg\|
 \bigg(\dsum_{j\in \Bbb Z}
 |T^*_{\Omega;b,u} f_j|^2\bigg)^{1/2}\bigg\|_{L^p}\le C\|b\|_{\ast}^u\|\Omega\|_{L\log^+L(\mathbf{S}^{n-1})}\bigg\|
 \bigg(\dsum_{j\in \Bbb Z}
 | f_j|^2\bigg)^{1/2}\bigg\|_{L^p}
$$
which were established in \cite{DR86} and \cite{CD}.
Then by  \eqref{S2kdm}  and Lemma \ref{Deltacom}, we
have for $1<p<\infty,$
\begin{align*}
\|S_{2,k}(\mathcal{T}_bf)\|_{L^p}&\le C\|b\|_{\ast}\|\Omega\|_{L\log^+L(\mathbf{S}^{n-1})}\|f\|_{L^p}.\end{align*}
This gives \eqref{S2kp}. Therefore, we complete the proof Proposition \ref{S2k}. \qed

\section{Proof of Theorem \ref{thm:N}}

We write
\begin{align*}
\Omega(x')&=[\Omega(x')-\frac1{\omega_{n-1}}\int_{\mathbf S^{n-1}}\Omega(y')d\sigma(y')]+\frac1{\omega_{n-1}}\int_{\mathbf S^{n-1}}\Omega(y')d\sigma(y')\\
&:=\Omega_0(x')+C(\Omega,n),
\end{align*}
where $\omega_{n-1}$ denotes the area of $\mathbf S^{n-1}$. Thus,
\begin{align*}
M_{\Omega,t;b}f(x)&=\frac1{t^n}\int_{|x-y|<t}\Omega_0(x-y)(b(x)-b(y))f(y)dy+C(\Omega,n)\frac1{t^n}\int_{|x-y|<t}(b(x)-b(y))f(y)dy\\&:=M_{\Omega_0,t;b}f(x)+C(\Omega,n)M_{t;b}f(x),
\end{align*}
where $\Omega_0$ satisfies the cancelation condition \eqref{can of O}. Denote the operator family $\{M_{\Omega_0,t;b}\}_{t>0}$ by $\mathcal M_{\Omega_0;b}$ and $\{M_{t;b}\}_{t>0}$ by $\mathcal M_b$.
By Corollary\ref{cor:4},  we get  for $1<p<\infty,$$$\|V_\rho\mathcal M_b(f)\|_{L^p}\le C\|b\|_{\ast}\|f\|_{L^p}.$$
 To prove Theorem \ref{thm:N}, it suffices to show
 \begin{equation}\label{NMO0}
 \|V_\rho\mathcal M_{\Omega_0;b}(f)\|_{L^p}\le C\|b\|_{\ast}\|\Omega\|_{L\log^+L(S^{n-1})}\|f\|_{L^p}.
   \end{equation}
   Similarly, the proof of  \eqref{NMO0} is  reduced to prove
\begin{equation}\label{NM0bp}
 \|V_\rho(\{M_{\Omega_0,2^k;b}f\}_{k\in \Bbb Z})\|_{L^p}\le C\|b\|_{\ast}\|\Omega\|_{L\log^+L(S^{n-1})}\|f\|_{L^p}.
\end{equation}
and
\begin{equation}\label{S2MO0bp}
\|S_2(\mathcal M_{\Omega_0;b}f)\|_{L^p}\le
C\|b\|_{\ast}\|\Omega\|_{L(\log^+L)^2(S^{n-1})}\|f\|_{L^p}.
\end{equation}
For \eqref{NM0bp}, the pointwise domination
$$
V_\rho(\{M_{\Omega_0,2^k;b}f\}_{k\in \Bbb Z})\leq (\sum_{k}|M_{\Omega_0,2^k;b}f|^2)^{1/2}
$$
for $1<p<\infty$,  which is a known result in \cite{HY}.

For \eqref{S2MO0bp},  observe that
$V_{2,j}(\mathcal{M}_{\Omega_0;b}f)$ is just the strong  $2$-variation function of the family
$\{M_{\Omega_0,2^jt;b} f\}_{t\in[1,2)}$, hence
\begin{align*}
S_{2}(\mathcal{M}_{\Omega_0;b}f)(x)&=\Big(\sum_{j\in\mathbb{Z}}|V_{2,j}(\mathcal{M}_{\Omega_0;b}f)(x)|^2\Big)^{\frac{1}{2}}.
\end{align*}
 Similar to the proof of \eqref{short variation singular}, we get that for $1<p<\infty,
 $\begin{align*}
\|S_{2}(\mathcal{M}_{\Omega_0;b}f)\|_{L^p}\le C\|b\|_{\ast}\|\Omega_0\|_{L(\log^+L)^2(S^{n-1})}\|f\|_{L^p}\leq C\|b\|_{\ast}\|\Omega\|_{L(\log^+L)^2(S^{n-1})}\|f\|_{L^p}.
\end{align*}
Therefore, \eqref{S2MO0bp} is proved.\qed

\bibliographystyle{amsplain}

\end{document}